\documentclass[11pt]{article}

\usepackage{amsmath,amssymb,enumerate,epsfig}

\newtheorem{proposition}{Proposition}
\newtheorem{theorem}{Theorem}
\newtheorem{lemma}{Lemma}

\newcommand{\byeproof}{\hfill $\square$}

\newcommand{\beginEL}[1]{\begin{equation}\label{#1}}
\newcommand{\closeE}{\end{equation}}

\newcommand{\bbR}{{\mathbb{R}}}

\newcommand{\U}{{\boldsymbol U}}
\newcommand{\V}{{\boldsymbol V}}
\renewcommand{\P}{{\boldsymbol P}}
\newcommand{\argmin}{{\mathrm{argmin}}}
\newcommand{\grad}{\nabla}
\newcommand{\pd}{\partial}
\newcommand{\bn}{{\boldsymbol n}}
\newcommand{\bm}{{\boldsymbol m}}
\newcommand{\bv}{{\boldsymbol v}}
\newcommand{\bone}{{\mathbf{1}}}
\newcommand{\be}{{\boldsymbol{e}}}
\newcommand{\bw}{{\boldsymbol w}}
\newcommand{\mean}[1]{\langle #1 \rangle}
\newcommand{\inner}[2]{\langle #1, #2 \rangle}

\newcommand{\eps}{\varepsilon}

\renewcommand{\H}{\mathcal{H}}
\newcommand{\OG}{{\Omega\setminus\Gamma}}
\newcommand{\Prob}{\mathrm{Prob}}
\newcommand{\E}{\mathrm{E}}
\newcommand{\sumK}{\sum_{i=1}^K}

\begin{document}

\setlength{\baselineskip}{13.2pt}

\title{A Stochastic-Variational Model for Soft Mumford-Shah Segmentation
       \thanks{This work has been partially supported by the NSF (USA) under grant number
            DMS-0202565. Email: jhshen@math.umn.edu; Tel: (USA) (612) 625-3570. }}

\author{Jianhong (Jackie) Shen \\
        {\small School of Mathematics,
        University of Minnesota,
        Minneapolis, MN 55455, USA} \\
        {\small Lotus Hill Institute for Computer Vision and Information Science,
         E'Zhou, Wuhan 436000, China}
        }

\date{}

\maketitle

\begin{abstract}
In contemporary image and vision analysis, stochastic approaches demonstrate great
flexibility in representing and modeling complex phenomena, while variational-PDE methods
gain enormous computational advantages over Monte-Carlo or other stochastic algorithms.
In combination, the two can lead to much more powerful novel models and efficient
algorithms. In the current work, we propose a stochastic-variational model for {\em soft}
(or fuzzy) Mumford-Shah segmentation of mixture image patterns. Unlike the classical {\em
hard} Mumford-Shah segmentation, the new model allows each pixel to belong to each image
pattern with some probability. We show that soft segmentation leads to hard segmentation,
and hence is more general. The modeling procedure, mathematical analysis, and
computational implementation of the new model are explored in detail, and numerical
examples of synthetic and natural images are presented.
\end{abstract}


{\small \tableofcontents}


\section{Introduction: Soft vs. Hard Segmentation}

Segmentation is the key step towards high-level vision modeling and analysis, including
object characterization, detection, and classification. There have been some recent
developments indicating that certain high-level visual tasks such as global scene
interpretation might be able to bypass segmentation~\cite{feifei99,feifei05}.
Nevertheless, segmentation still remains perhaps {\em the} most important and inspiring
task to date in low- or middle-level vision analysis and image processing.

The segmentation problem can be formulated as follows. Given an image $I$ on a
2-dimensional (2D) domain $\Omega$ (assumed to be bounded, smooth, and open), one seeks
out a {\em closed} ``edge set" $\Gamma$, and all the connected components $\Omega_1,
\dots, \Omega_K$ of $\OG$, such that by certain suitable visual measure (e.g., textural
or photometric), the image $I$ is discontinuous along $\Gamma$ while smooth or
homogeneous on each segment $\Omega_i$. Each image patch $I_i= I\big|_{\Omega_i}$ is also
called a {\em pattern}, and $\Omega_i$ its {\em support}.

We shall call this most common practice ``hard" segmentation. A hard segmentation
partitions the image domain $\Omega$ along a {\em definitive} edge set $\Gamma$, and
outputs {\em non-overlapping} pattern supports $\Omega_1, \dots, \Omega_K$.

The present work introduces the notion of ``soft" segmentation. Mathematically, a hard
segmentation amounts to the partition of the unit using indicator functions:
 \[
 1_\Omega(x) = \sum_{i=1}^K 1_{\Omega_i} (x), \qquad a.e. \quad x=(x_1, x_2) \in \Omega.
 \]
A soft segmentation seeks out instead a softer partition of the unit:
 \beginEL{E1:Spou}
 1_\Omega(x) = \sum_{i=1}^K  p_i(x),
 \closeE
where $p_i$'s are continuous or smoother functions. Formally, each $p_i$ could be
considered as the mollified version of $1_{\Omega_i}(x)$.

In the stochastic literature of image analysis and modeling, the above notion of soft
segmentation is closely connected to {\em mixture image models}
(e.g.,~\cite{JepsonBlack}). Suppose a given image $I$ is composed from $K$ unknown
patterns:
 \[ \omega=1,  \quad \omega=2, \quad \dots, \quad \omega=K, \]
where $\omega$ denotes the pattern label variable. At each pixel $x \in \Omega$,
$\omega(x) \in \{1, \dots, K\}$ becomes a random variable. Then the $p_i$'s in
(\ref{E1:Spou}) carry the natural stochastic interpretation:
 \[
 p_i(x) = \Prob(\omega(x)=i), \qquad i=1:K.
 \]
For this reason, each $p_i$ shall be called the {\em ownership} of pattern $i$. Instead
of the repulsive ownership in a hard segmentation, a soft one allows each pattern to
``own" a pixel with some likelihood.

Soft segmentation is more general since it can lead to natural hard segmentation under
the {\em maximum likelihood} (ML) principle. Given a soft segmentation $\{p_i(x):
i=1:K\}$, one can define for each pixel $x \in \Omega$ its unique owner $\omega_\ast(x)$
by:
 \beginEL{E1:s2hA}
 \omega_\ast(x) = \mathrm{argmax}_{\omega \in 1:K} \; p_\omega(x),
 \closeE
and if the maxima are non-unique, accept the largest index from the argmax pool. The
segments are then defined by
 \beginEL{E1:s2hB}
 \Omega_i = \omega_\ast^{-1}(i)=\{ x \in \Omega \mid \omega_\ast(x)=i \}, \quad i=1:K,
 \closeE
which leads to a natural hard segmentation. (\ref{E1:s2hA}) and (\ref{E1:s2hB}) are
called the {\em hardening formulae}.

Soft segmentation has been motivated by practical analysis of natural images. Patterns in
natural scenes often do not have clear-cut boundaries. In Figure~\ref{F1:soft}, for
example, there does not seem to exist a ``hard" boundary between the grass and sand
areas. If one draws an oriented line as shown in the figure, it makes more sense to state
that along the arrow the pattern transits from being ``more" sand like to being ``more"
grass like. Such consideration favors the following stochastic view that along the arrow,
the ownership
 \[
 \Prob(\omega(x)=\mathrm{grass}) \;\; \mbox{increases, \;\; while \;\;}
 \Prob(\omega(x)=\mbox{sand}) \;\; \mbox{decreases}.
 \]

{
\begin{figure}[ht]
 \centering{
 \epsfig{file=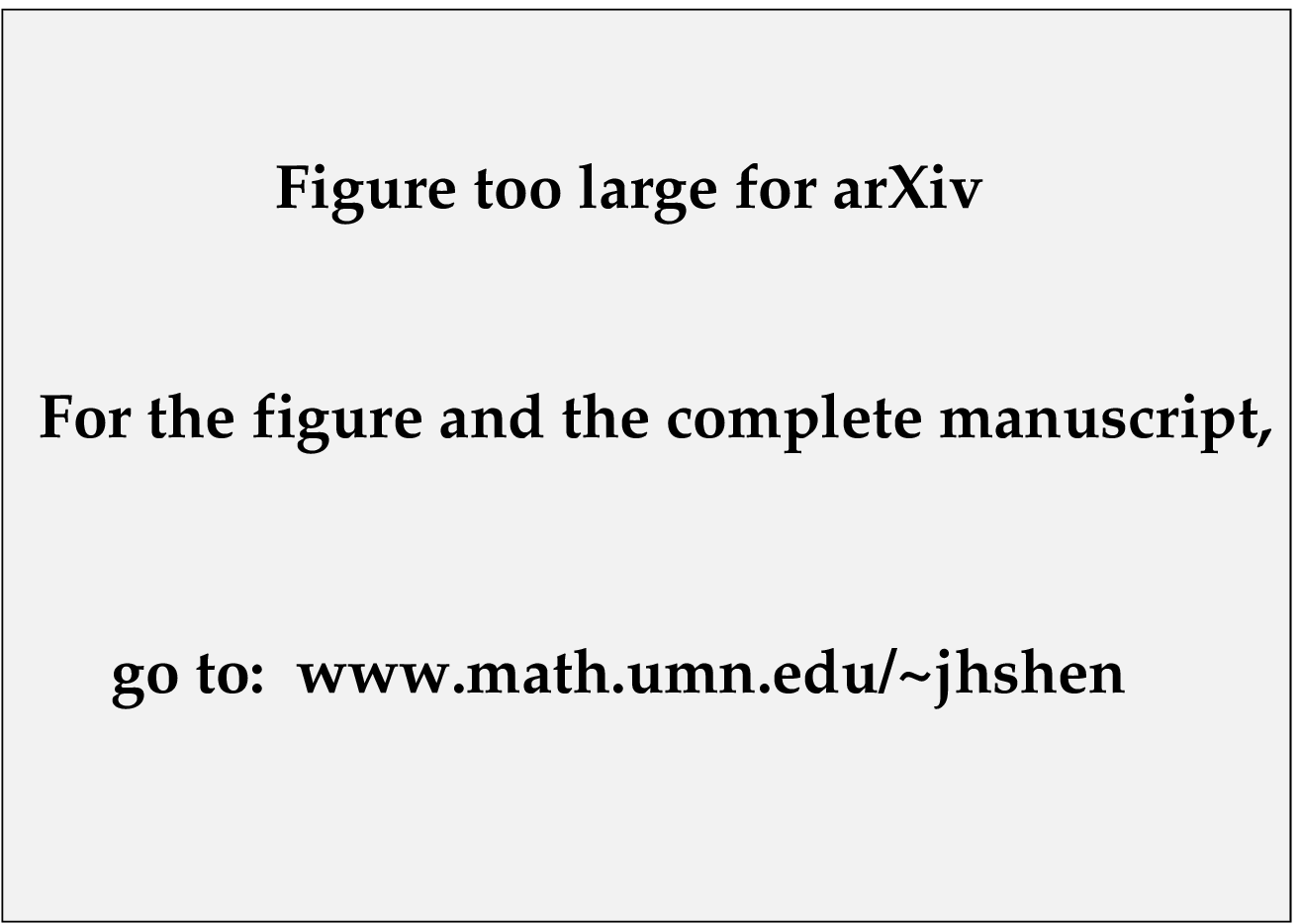, height=2in}
 \caption{Natural images often do not have clear-cut ``hard" boundaries between
          different patterns. Along the arrow, for example, one only observes that
          the sand pattern {\em gradually} becomes a grass pattern.
          Such a ``soft" view is the stochastic view on the segmentation problem.}
\label{F1:soft} }
\end{figure}
}

In the present work, we propose a new stochastic-variational soft-segmentation model for
the following celebrated Mumford-Shah model~\cite{morsol,mumsha}:
 \beginEL{E1:HMS}
 \min_{\Gamma, u} E[u, \Gamma \mid I]
  = \H^{1}(\Gamma) + \alpha \int_\OG |\grad u|^2 + \lambda \int_\Omega (u -I)^2,
 \closeE
where $\H^1$ stands for the 1D Hausdorff measure, which is simply the length when
$\Gamma$ is regular enough. {\em For notational conciseness, the default area-element
symbol $dx=dx_1dx_2$ will be omitted in most integral formulae.}

As stated in the abstract, the stochastic softness induces more flexibility and
universality in modeling, while the variational-PDE approach facilitates rigorous
mathematical analysis as well as more efficient computational implementations compared
with purely stochastic approaches including, e.g., the Monte-Carlo method or Gibbs'
sampling~\cite{bre,gemgem,tu_zhu,zhuyui}.

The paper has been organized as follows. Section~2 builds up the soft Mumford-Shah (SMS)
model under the Bayesian rationale and the MAP estimator~\cite{kniric,mum_bayes}, which
are the formal stochastic foundations of the present model. In Section~3, the prior
energy on the ownerships $p_i$'s is developed based on the celebrated work of Modica and
Mortola~\cite{modmor} on phase-field modeling and $\Gamma$-convergence approximation in
material sciences and phase transitions. In Section~4, we analyze the main mathematical
properties of the proposed SMS model, including the admissible space, hidden symmetry and
symmetry breaking via weak supervision, and the existence theorems. In Section~5, we then
derive the system of Euler-Lagrange equations of the SMS model for which the role of the
probability simplex constraint is discussed in detail. Section~5 also introduces the
alternating-minimization algorithm to compute the Euler-Lagrange equations. Finally, the
numerical performance of the SMS model is demonstrated in Section~6 via both synthetic
and natural test images that are sufficiently representative and generic.

Throughout the manuscript, the notation $E[X, Y \mid Z]$ in the deterministic setting
always denotes a quantity (often a functional) $E$ that depends on $X, Y,$ and $Z$ but
with $Z$ given or fixed. Similarly, $E[X \mid Y, Z]$ still denotes $E[X, Y \mid Z]$
modulo some additive quantity $g[Y, Z]$ that is often unimportant as far as the
optimization on $X$ (given $Y$ and $Z$) is concerned. These notations therefore have been
closely inspired by {\em conditional probabilities} in the stochastic setting.


\section{Bayesian Rationale to the New Model and Gaussian Mixture}

\subsection{Bayesian Rationale}

Segmentation can be done in some feature spaces such as gradient-like highpass features
or Gabor features (e.g.,~\cite{cvs3,zhuwumum,zhuyui}). The Mumford-Shah model easily
extends to such general features (e.g.,~\cite{cvs3}), even though it was originally
formulated only for intensity fields. For maximal clarity in exposing the core ideas of
the current work, we shall also focus only on the latter, while leaving as canonical
exercises to adapt the new model for any given feature distribution.

Let $K$ be the total number of intended patterns. As in~\cite{tu_zhu,zhuyui}, $K$ could
also be treated as an unknown to be optimally estimated, which however does not add much
to the most significant contribution (i.e., the modeling and computation of the
``softening" procedure) of the present work.

Given an image input $I=I(x)$ on a bounded, regular, and open domain $\Omega$, the
primary goal of soft segmentation is to compute the ownerships:
 \[ p_1(x), p_2(x), \dots, p_K(x). \]
Define $\P(x)=(p_1(x), p_2(x), \dots, p_K(x))$, and
 \[ \Delta_{K-1} = \mbox{convex hull of} \; \vec{e}_1, \dots, \vec{e}_K, \]
where the $(\vec{e}_i \mid i=1:K)$ denotes the canonical Cartesian basis of $\bbR^K$.
$\Delta_{K-1}$ is often called the {\em canonical} $(K-1)$-simplex, or the {\em
probability}-simplex in $\bbR^K$. Then
 \[
  \P: \Omega \to \Delta_{K-1}, \qquad x \to \P(x),
 \]
meaning that the total ownerships always add up to 100\% at any pixel $x\in \Omega$.

Associated with each pattern label $\omega=i$ is a smooth function $u_i(x) \in
H^1(\Omega)$, similar to the original Mumford-Shah model. Here the Sobolev space
$H^1(\Omega)$ is defined by~\cite{shen_MSS}
 \[
 H^1(\Omega) = \{ u \in L^2(\Omega) \mid \grad u \in L^2(\Omega, \bbR^2) \}.
 \]
Define $\U(x)=(u_1(x), u_2(x), \dots, u_K(x))$. Then the goal of soft segmentation is to
estimate the optimal vectorial pair of ownerships and patterns given an image $I$:
 \[
 (\P_\ast, \U_\ast) = \mathrm{argmax}_{(\P,\U)} \Prob (\P, \U \mid I).
 \]

By the Bayesian formula~\cite{kniric,mum_bayes}, the posterior given $I$ is expressible
via
 \[
 \Prob(\P, \U \mid I) = \Prob(I \mid \P, \U) \Prob(\P) \Prob(\U) / \Prob(I),
 \]
assuming that the mixture patterns $\U$ and the mixture rules $\P$ are {\em independent}.
We shall call the first term a ``mixture generation" model, since it reveals how the
image data should look like given the information of the patterns and their ownerships.

By taking the logarithmic likelihood $E[\cdot]=-\log \Prob(\cdot)$, or formally Gibbs'
energy in statistical mechanics~\cite{cha_IMSM,gib}, one attains the soft segmentation
model in its ``energy" form:
 \beginEL{E2:V-SMS}
 \min_{(\P, \U)} E[\P, \U \mid I] =
   E[I \mid \P, \U] + E[\P] + E[\U],
 \closeE
modulo an insignificant additive constant.

Assuming that all the pattern channels are independent from each other, one has
 \[ E[\U] =E[u_1, \dots, u_K] = \sum_{i=1}^K E[u_i \mid i]. \]
For Sobolev-regular patterns, one may impose the homogeneous Sobolev energies:
 \beginEL{E2:sobolev}
 E[u_i \mid i] = E[u_i]= \alpha \int_\Omega |\grad u_i|^2  , \qquad i=1:K,
 \closeE
for some scalar weight $\alpha$ that models the visual sensitivity to intensity
roughness. Unlike the original Mumford-Shah model, the energy for each channel has been
defined on the entire image domain $\Omega$ instead of on each ``hard-cut" patch
$\Omega_i$. Thus the energy form~(\ref{E2:sobolev}) must carry out {\em extrapolation}
for practical applications. Long-range extrapolations are, however, often unimportant
after being weighed down by their negligible ownerships $p_i$'s.

\subsection{Gaussian Mixture with Smooth Mean Fields}

In this subsection we discuss the mixture generation model $\Prob(I \mid \P,\U)$ or $E[I
\mid \P, \U]$.

Assume that the patterns are all Gaussian's with mean fields $u_1, u_2, \dots, u_K$. For
simplicity also assume that they share the same variance $\sigma^2$ (which readily
generalizes to the more general case with variations). Then at any given pixel $x \in
\Omega$,
 \[ (I \mid \omega(x)=i) \sim N(u_i(x), \sigma^2), \qquad i=1:K. \]
Define the Gaussian probability density function (p.d.f)
 \[
 g(I \mid m, \sigma) = \frac 1 {\sqrt{2\pi}\sigma} \exp\left(-\frac {(I-m)^2} {2 \sigma^2} \right).
 \]
The the p.d.f of the mixture image $I$ at any pixel $x$ is given by
 \[
  \Prob \left( I(x) \mid \P(x), \U(x) \right)
            = \sum_{i=1}^K \Prob(I \mid \omega(x)=i) \; \Prob(\omega(x)=i)
            = \sum_{i=1}^K g(I \mid u_i(x), \sigma) p_i(x).
 \]
Thus {\em ideally} the ``energy" for the mixture generation model should be given by
 \beginEL{E2:EmixGau}
 E[I \mid \P, \U] = -  \mu \int_\Omega \log
 \left( \sum_{i=1}^K g(I \mid u_i(x), \sigma) p_i(x) \right), \quad
 \mbox{for some scalar $\mu>0$,}
 \closeE
provided that, given two fields $\P$ and $\U$ on $\Omega$, for any two distinct pixels
$x$ and $y$,
 \[
  (I(x) \mid \P, \U) \quad \mbox{is independent of} \qquad (I(y) \mid \P, \U).
 \]

In the current work, we shall adopt a reduced form of the complex
formula~(\ref{E2:EmixGau}), which is simpler and easier to manage both in theory and for
computation. Assume that each soft ownership $p_i(x)$ is closer to a hard one $p_i(x)
\simeq 1_{\Omega_i}(x)$ for $i=1:K$. Then
 \[
 \begin{split}
 - \log \left( \sum_{i=1}^K g(I \mid u_i(x), \sigma) p_i(x) \right)
 & \simeq - \log \left( \sum_{i=1}^K g(I \mid u_i(x), \sigma) 1_{\Omega_i} (x) \right) \\
 & = - \sum_{i=1}^K \log g(I \mid u_i(x), \sigma) 1_{\Omega_i}(x)  \qquad (a.\, e.) \\
 & \simeq - \sum_{i=1}^K \log g(I \mid u_i(x), \sigma) p_i(x) \\
 & = \frac 1 {2\sigma^2} \sum_{i=1}^K (I - u_i(x))^2 p_i(x) + const., \\
 \end{split}
 \]
where the additive constant only depends on $\sigma$ and $K$. This suggests the following
convenient energy form for the mixture generation model:
 \beginEL{E2:mixgen}
 E[I \mid \P, \U] = \lambda \int_\Omega
    \left( \sum_{i=1}^K (I - u_i(x))^2 p_i(x) \right)  ,
 \closeE
which amounts to a weighted least-square energy~\cite{str}. The weight $\lambda$ reflects
visual sensitivity to synthesis errors.

In combination of (\ref{E2:V-SMS}), (\ref{E2:sobolev}), and (\ref{E2:mixgen}), the new
soft segmentation model takes the form of minimizing
 \beginEL{E2:combined}
  E[\P,\U \mid I]=
  \lambda \sumK \int_\Omega (I -u_i(x))^2 p_i(x)   + \alpha \sumK \int_\Omega |\grad
  u_i|^2   + E[\P].
 \closeE
Notice that here the ownership distribution $\P$ ``softens" the ``hard" segmentation
boundary $\Gamma$ in the original Mumford-Shah model~(\ref{E1:HMS}). To complete the
modeling process, it suffices to properly define the prior or regularity energy $E[\P]$,
which is the main task of the next section.

\section{Modica-Mortola's Phase-Field Model for Ownership Energy}

To generalize but not to deviate too far from classical hard segmentation, it is natural
to impose the following two constraints:
 \begin{enumerate}[(a)]
 \item each pattern ownership $p_i(x)$ has almost only two phases: on (corresponding to
 $p_i=1$) and off (to $p_i=0$), and the transition band in between is narrow; \quad and

 \item the soft boundaries, or equivalently the transition bands, are {\em regular},
 instead of being zigzag.
 \end{enumerate}
In combination, one imposes the following Modica-Mortola type of energy with a
double-well potential~\cite{modmor}: $p_i \in H^1(\Omega)$,
 \beginEL{E3:MMenergy}
 E_\eps[p_i]=\int_\Omega \left( 9\eps|\grad p_i|^2 +{(p_i(1-p_i))^2\over\eps}\right) ,
 \quad i=1:K.
 \closeE
Here $\eps \ll 1$ controls the transition bandwidth. Since $\eps \ll 1$, the second term
necessarily demands $p_i \simeq 0$ or $1$ to lower the energy, which well resonates with
the expectation in (a). The first term, weighted by the small parameter $\eps$, amounts
to a regularity condition on each $p_i$, which meets the requirement in (b).

Energies in the form of (\ref{E3:MMenergy}) are very common in material sciences,
including the theories of liquid crystals and phase transitions~\cite{eri,ginlan}.
Mathematically they have been well studied in the framework of
$\Gamma$-convergence~\cite{DalMasoBook}, which we now give a brief introduction in the
present context. We also refer the reader to the works of Ambrosio and
Tortorelli~\cite{ambtor1,ambtor2} on the $\Gamma$-convergence approximation to the
classical Mumford-Shah segmentation model.

Recall that for any $q(x) \in L^1(\Omega)$, its {\em total variation} as a Radon measure
is defined by~\cite{chashe_book,giu}
 \[
 \mathrm{TV}[q]=\int_\Omega |Dq| =\sup_{\boldsymbol{g} \in C^1_0(\Omega,B^2)}
                \inner{q}{\nabla \cdot \boldsymbol{g}},
 \]
where $B^2$ stands for the unit disk centered at the origin in $\bbR^2$. (The TV measure
was first introduced into image processing by Rudin, Osher, and Fatemi~\cite{rudoshfat}.)
Define for any $q \in L^1(\Omega)$,
\[
 E_0[q] =
          \begin{cases}
            \mathrm{TV}[q], & \qquad \mbox{if $q=0$ or $1$, $a.e.$ on $\Omega$,}  \\
            \infty,         & \qquad \mbox{otherwise.}
          \end{cases}
\]
As a result, a finite energy $E_0[q]$ necessarily implies that $q$ has two phases only,
and $E_0[q]=\mathrm{TV}[q]=\mathrm{Per}(q^{-1}(1))$ is the perimeter of the support
region $V=q^{-1}(1)$.

Further define
 \[ L^1_{[0,1]}(\Omega) =\{ q \in L^1(\Omega) \mid q(x) \in [0, 1], \; \forall \; x \in \Omega \} \]
to be a subspace of $L^1(\Omega)$ (as a metric space). Then Modica and Mortola's well
known results in~\cite{modmor} readily leads to the following theorem.

\begin{theorem}[$\Gamma$-Convergence Approximation of a Two-Phase $\mathrm{TV}$] For any
$q \in L^1_{[0,1]}(\Omega) \setminus H^1(\Omega)$, we extend the definition of
$E_\eps[\cdot]$ in (\ref{E3:MMenergy}) by defining $E_\eps[q]=+\infty$. Then
 \[
    E_\eps \to E_0 \qquad \mbox{in the sense of $\Gamma$-convergence in
    the metric space $L^1_{[0,1]}(\Omega)$.}
 \]
That is
\begin{enumerate}[{\rm (i)}]
 \item for any $q_\eps \to q$ in $L^1_{[0,1]}(\Omega)$ as $\eps \to 0$,
        \[ \liminf_{\eps \to 0} E_\eps[q_\eps] \ge E_0[q]; \qquad \mbox{and} \]
 \item for any $q \in L^1_{[0,1]}(\Omega)$, there exists some sequence
       $(q_\eps^\ast \mid \eps)$, such that $q^\ast_\eps \to q$ as $\eps \to 0$, and
        \[ \lim_{\eps \to 0} E_\eps[q^\ast_\eps] =E_0[q]. \]
\end{enumerate}
\end{theorem}

We refer the reader to Modica and Mortola~\cite{modmor} for a proof (with some necessary
modification). Here we only point out that the ``tight" sequence $(q^\ast_\eps \mid
\eps)$ in (ii) can be constructed using a smooth sigmoid transition across the hard
boundary of a given two-phase function $q$. Recall as in the theory of neural
networks~\cite{neuralnetBook} that a {\em sigmoid} transition between 0 and 1 is achieved
by
 \[ \sigma(t) = \frac {1}{1+ e^{-t}}, \qquad -\infty < t < \infty. \]
The scaling parameter $\eps$ participates in the transition by the form of
$\sigma(t/(3\eps))$. In particular, $\eps$ indeed corresponds to the width of the
transition band when $t$ is a distance function.

This theorem reveals the close connection of the particular choice of $E_\eps[p_i]$ in
(\ref{E3:MMenergy}) to the original Mumford-Shah model.

\begin{proposition}
Suppose $p_\eps$'s ``optimally" (i.e., by the above sigmoidal transition) converge to a
given 2-phase pattern $1_V(x)$ with a regular hard boundary $\Gamma=\pd V$. Then,
 \[ E_\eps[p_\eps] \to \mathrm{length}(\Gamma) = \int_\Omega \left|D 1_V(x)\right|. \]
\end{proposition}

Similar results have appeared in the earlier influential works of Ambrosio and
Tortorelli~\cite{ambtor1,ambtor2} on the $\Gamma$-convergence approximation to the
Mumford-Shah model. The technique has also been extensively applied in image computation
and modeling~\cite{eseshe,mar92,mardoz,shen_Gampcms,she_mirror,shejun_wmsbe} to overcome
the difficulty in representing and computing the free boundary $\Gamma$.

To summarize this section, we propose the following energy model for the ownership
distribution $\P(x)=(p_1(x), p_2(x), \dots, p_K(x))$:
 \beginEL{E3:KmmE}
 E_\eps[\P]=\sumK E_\eps[p_i]
    = \sumK \int_\Omega \left( 9\eps|\grad p_i|^2 + {(p_i(1-p_i))^2 \over \eps} \right).
 \closeE
One, however, must realize that different ownerships are {\em not} decoupled by this
energy though it has appeared so. The energy $E_\eps[\P]$ must be coupled with the
constraint of the probability-simplex:
 \[
     \P: \Omega \to \Delta_{K-1}, \qquad \mbox{or} \qquad \sumK p_i(x)\equiv 1, \quad
        p_i \ge 0, \quad \forall \; x \in \Omega.
 \]
In particular, for small $\eps$, although (\ref{E3:KmmE}) implies that each ownership
$p_i$ tends to polarize to $0$ or $1$ independently, they have to cooperate with each
other under the above simplex constraint to optimally share the ownerships.


\section{Soft Mumford-Shah (SMS) Segmentation}

\subsection{The Model and Admission Space}

Combining the preceding two sections, we have developed the complete formula for soft
Mumford-Shah segmentation with $K$-patterns:
 \beginEL{E4:SMSa}
 \min_{\P,\U} E[\P, \U \mid I ]
 =   \lambda \sumK \int_\Omega (u_i -I)^2 p_i
   + \alpha \sumK \int_\Omega |\grad u_i|^2
   + \sumK \int_\Omega \left( 9\eps|\grad p_i|^2 + {(p_i(1-p_i))^2 \over \eps} \right)
    ,
 \closeE
with the constraint that
 \[
 \P: \Omega \to \Delta_{K-1}, \qquad \mbox{the probability $(K-1)$-simplex},
 \]
i.e., $p_i \ge 0, i=1:K$, and $\sumK p_i =1$. As discussed previously, it is this simplex
constraint that induces coupling among different channels into the seemingly decoupled
model~(\ref{E4:SMSa}).

Besides the simplex constraint, the last term in the energy (\ref{E4:SMSa}) requires $p_i
\in H^1(\Omega)$ for $i=1:K$. Similarly, the second term requires each pattern $u_i \in
H^1(\Omega)$. Then with the assumption that
 \[     \mbox{``the given image $I \in L^2(\Omega)$,"}   \]
$E[\P, \U \mid I]$ is well defined and finite for any admissible patterns $\U$ and
pattern ownership distribution $\P$:
 \beginEL{E4:admK}
 \mathrm{adm}_K =
 \{
 (\P, \U) \mid p_i, u_i \in H^1(\Omega), i=1:K; \quad \P: \Omega \to \Delta_{K-1}
 \}.
 \closeE

\subsection{Breaking the Hidden Symmetry via Weak Supervision}

Let $S_K$ denote the permutation group of $\{ 1, \dots, K\}$. Each permutation $\sigma
\in S_K$ is a 1-to-1 map:
 \[ \sigma: \{1, \dots, K\} \to \{1, \dots, K\}, \]
so that $(\sigma(1), \dots, \sigma(K))$ is a re-arrangement of $\{1, \dots, K\}$. For any
$K$-tuple $\vec{F}=(f_1, \dots, f_K)$, one defines
 \[
 \vec{F}_\sigma = (f_{\sigma(1)}, f_{\sigma(2)}, \dots, f_{\sigma(K)}).
 \]

\begin{theorem}[Hidden Symmetry of SMS]
 For any $\sigma \in S_K$,
 \[ E[\P_\sigma, \U_\sigma \mid I] = E[\P, \U \mid I]. \]
 In particular, suppose
 \[
 (\P^\ast, \U^\ast) = \argmin_{(\P, \U) \in \mathrm{adm}_K} E[\P, \U \mid I]
 \]
 is an optimal pair. Then for any $\sigma \in S_K$, $(\P^\ast_\sigma, \U^\ast_\sigma)$
 is a minimizer as well.
\end{theorem}

The proof is straightforward and thus omitted. Such symmetry not only worsens the
non-uniqueness of the non-convex energy functional in (\ref{E4:SMSa}), but also
potentially jitters intermediate solutions in iterative computational schemes (i.e.,
hysterical transitions in the admissible space).

To break the permutation symmetry, we turn to a weak supervision scheme in which a user
specifies $K$ distinct domain patches:
 \[ Q_1, Q_2, \dots, Q_K, \]
and imposes the symmetry-breaking conditions:
 \beginEL{E4:supervision}
   p_i \big|_{Q_j} = \delta_{ij}, \qquad i, j=1:K,
 \closeE
where $\delta_{ij}$ denotes Kronecker's delta. That is, a user requires each given patch
$Q_i$ to be a ``pure" pattern exclusively labelled by $i$.  Computationally this weak
supervision process can be automated based on multiscale patch statistics as in the
contemporary works on scene recognition~\cite{feifei05,feifei99}, or more generally, the
learning theory~\cite{cuckersmale01,smazho}.

{
\begin{figure}[ht]
 \centering{
 \epsfig{file=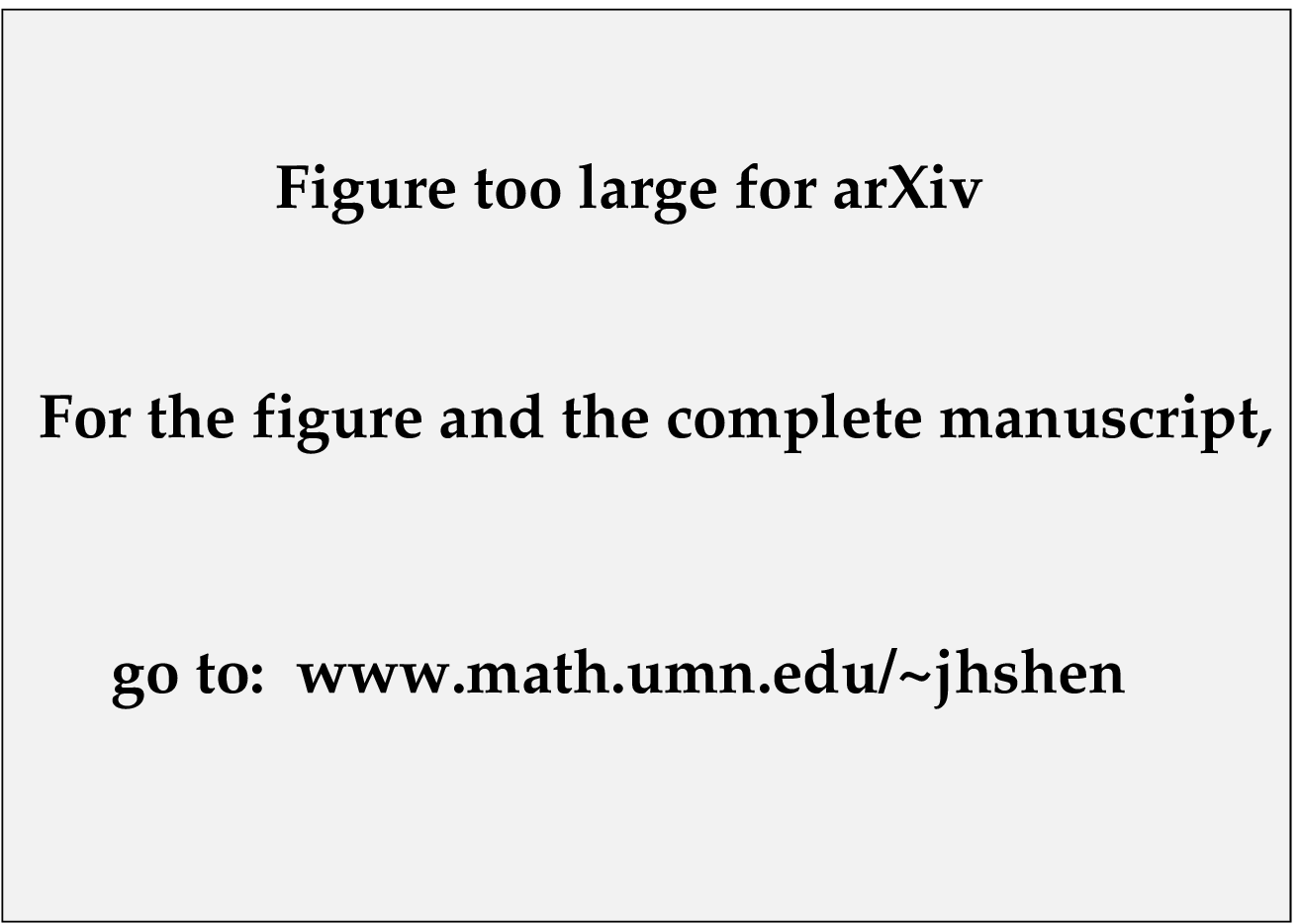, height=2in}
 \caption{Examples of a 3-phase supervision and a 4-phase supervision to break the
          symmetry in the model. Such weak supervision can also be automated based
          on multiscale patch statistics~\cite{feifei99,feifei05}. }
\label{F4:superbeach} }
\end{figure}
}

\subsection{Existence Theorems for Non-Supervision and Supervision}

In this subsection, we establish the existence theorems for the soft Mumford-Shah
segmentation model~(\ref{E4:SMSa}) without or with the supervision
(\ref{E4:supervision}).

\begin{theorem}[Existence Theorem for Unsupervised SMS]     \label{T4:unsuper}
Suppose $I \in L^2(\Omega)$. Then for any positive modeling parameters $(\lambda, \alpha,
\eps)$, a minimizer to the unsupervised soft Mumford-Shah model~(\ref{E4:SMSa}) must
exist.
\end{theorem}

We will need the following lemma for the proof.

\begin{lemma} \label{L4:meanbound}
Let $(f_n \mid n)$ be a sequence of functions in $L^2(\Omega)$, and $(p^n \mid n)$ a
sequence of nonnegative measurable functions on $\Omega$ and valued in $[0, 1]$.  Suppose
 \begin{enumerate}[{\rm (i)}]
 \item $p^n \to p^\ast$, a.e. on $\Omega$, and $\int_\Omega p^\ast > 0$; and
 \item $\int_\Omega f_n^2 p^n \le A$ for some $A>0$ and $n=1:\infty$.
 \end{enumerate}
Then there exists some function $\rho \in L^2(\Omega)$, such that
 \begin{enumerate}[{\rm (a)}]
 \item $\rho \ge 0$ and $\int_\Omega \rho =1$, and
 \item for some fixed $B>0$, $\left| \int_\Omega f_n \rho \right| \le B$ for $n=1:\infty$.
 \end{enumerate}
\end{lemma}

\begin{proof} Denote the Lebesgue measure of a measurable set $W$ by $|W|$.
Since $p^\ast \ge 0$ and $\int_\Omega p^\ast >
0$, there must exist some $c>0$, such that
 \[ V= \{ x\in \Omega \mid p^\ast > 2 c\} \quad \mbox{has a finite but positive measure.} \]
On the other hand, by Egorov's theorem~\cite{fol} on a.e. convergence, there must exist a
subset $W \subset V$, such that
 \begin{enumerate}[(a')]
    \item $|V -W| \le \frac {|V|} 2$, and hence $|W|>0$, and
    \item $p^n \to p^\ast$ uniformly on $W$.
 \end{enumerate}
In particular, there exists some $N$, such that for any $n>N$, $p^n > c$ on $W$. Define
 \[ \rho(x) = { {1_W(x)} \over |W| } \in L^2(\Omega). \]
Then $\int_\Omega \rho =1$, and for any $n > N$,
 \[
 \int_\Omega f_n^2 \rho= {1 \over c|W| } \int_W f_n^2 c \le {1 \over c|W|} \int_\Omega
 f_n^2 p^n \le {A \over c|W|}.
 \]
 Thus by the Schwarz inequality (or $\E[X]^2 \le \E[X^2]$ in probability theory),
  \[
         |\int_\Omega f_n \rho | \le
 \left(  \int_\Omega f_n^2 \rho \right)^{1/2} \le \left(A \over c|W|\right)^{1/2}, \quad n >N
  \]
The lemma holds if one defines
 \( B= \max \left( \left(A \over c|W|\right)^{1/2}, \left|\int_\Omega f_1 \rho \right|,
            \dots, \left|\int_\Omega f_N \rho \right|\right).
 \)\byeproof
\end{proof}

We are ready to prove Theorem~\ref{T4:unsuper}.

\begin{proof}
Take the special pattern distribution:
  \[ u_i \equiv 0, \quad i=1:K; \quad p_1\equiv 1, \quad p_j\equiv 0, \quad j=2:K. \]
Then
 \[ E[\P, \U \mid I] = \lambda \int_\Omega I^2   < \infty. \]
Thus the infimum of the energy must be finite. Let $(\P^n, \U^n \mid n) \subseteq
\mathrm{adm}_K$ (see (\ref{E4:admK})) be a minimizing sequence for the soft Mumford-Shah
energy~(\ref{E4:SMSa}).

Due to the third term in the energy and the simplex constraint, for each channel $i$,
$(p^n_i \mid n)$ must be bounded in $H^1(\Omega)$. By the $L^2$-weak compactness, there
must exist some $\P^\ast \in L^2(\Omega, \bbR^K)$, and a subsequence of $(\P^n \mid n)$,
which after relabelling shall still be denoted by $(\P^n \mid n)$ for convenience, such
that
 \[ \P^n \to \P^\ast \quad \mbox{in $L^2(\Omega, \bbR^K)$}, \quad  n \to \infty. \]
Then by the $L^2$ lower semi-continuity of Sobolev measures,
 \beginEL{E4:lowerA}
 9 \eps \int_\Omega |\grad p^\ast_i|^2
   \le \liminf_{n\to \infty} \; \; 9\eps \int_\Omega |\grad p^n_i|^2, \qquad i=1:K.
 \closeE
Furthermore, with possibly another round of subsequence refinement, one can assume
 \[
 \P^n (x) \to \P^\ast (x), \quad a.e. \quad x \in \Omega, \qquad n \to \infty.
 \]
Since the probability simplex $\Delta_{k-1}$ is closed and $\P^n(x) \in \Delta_{K-1}$,
one concludes that
 \[
 \P^\ast(x) \in \Delta_{K-1}, \qquad a.e. \quad x \in \Omega.
 \]
And by Fatou's Lemma~\cite{fol,lielos}, one has
 \beginEL{E4:lowerC}
 \int_\Omega {(p^\ast_i (1-p^\ast_i))^2 \over \eps } \le
 \liminf_{n \to \infty} \int_\Omega {(p^n_i (1-p^n_i))^2 \over \eps }, \qquad i=1:K.
 \closeE
(In fact, the equality holds by {\em Lebesgue's Dominated Convergence}~\cite{lielos}.)

After the above subsequence selection on $\P^n$'s, one naturally has an associated
subsequence of $(\U^n \mid n)$, which for convenience is still denoted by $(\U^n \mid n)$
after relabelling. For each {\em specific} channel $i$, we then consider two scenarios
separately.

Suppose $p^\ast_i(x) \equiv 0, a.e. \; x \in \Omega$. We then define for that channel
 \beginEL{E4:lowerD}
  u^\ast_i(x) \equiv 0, \qquad x \in \Omega.
 \closeE
Such a channel is called a ``dumb" channel.

Otherwise, one must have $\int_\Omega p^\ast_i > 0$, and from the first term
in~(\ref{E4:SMSa}),
 \[
  \int_\Omega (u^n_i -I)^2 p^n_i  \le const., \quad n=1:\infty.
 \]
Since $\int_\Omega I^2 p^n_i \le \int_\Omega I^2$, by the triangle inequality,
 \[
 \int_\Omega (u^n_i)^2 p^n_i \le const., \qquad n=1:\infty,
 \]
where the constant only depends on $I$ and the model parameters. Then by
Lemma~\ref{L4:meanbound}, there exists some $\rho_i(x) \ge 0$, with $\int_\Omega \rho_i
=1$, some constant $B_i >0$ such that
 \[
 \left| \int_\Omega u^n_i \rho_i \right| \le B_i, \qquad n=1:\infty.
 \]
On the other hand, by the second term in the energy~(\ref{E4:SMSa}),
 \[
 \int_\Omega |\grad u^n_i|^2  \le C_i=C_i(I, \lambda, \alpha, \eps), \quad n=1:\infty,
 \]
for some constant $C_i$ independent of $n$. Then by the generalized Pointcar\'{e}
inequality~\cite{eva98,lielos} on $\Omega$,
 \[
 \|w - \inner{w}{\rho_i}\|_{L^2} \le A_i \|\grad w\|_{L^2},
 \]
where $A_i=A_i(\rho_i, \Omega)$ is independent of $w \in H^1(\Omega)$, one concludes that
 \[ \| u^n_i \|_{L^2} \le D_i=D_i(A_i, B_i, C_i), \qquad n=1:\infty, \]
for some constant $D_i$. As a result, $(u^n_i \mid n)$ must be bounded in $H^1(\Omega)$.
By the $L^2$-weak compactness of bounded $H^1$-sequences, there is a subsequence of
$(u^n_i \mid n)$, for convenience still denoted by $(u^n_i \mid n)$ after relabelling,
such that
 \[ u^n_i \to u^\ast_i \in L^2(\Omega), \qquad n \to \infty, \]
converging in the sense of both $L^2$ and almost everywhere. Then by the lower
semi-continuity,
 \[
 \int_\Omega |\grad u^\ast_i|^2 \le \liminf_{n\to \infty} |\grad u^n_i|^2.
 \]
Finally, since $u^n_i(x) \to u^\ast_i(x)$ and $p^n_i(x) \to p^\ast_i(x)$, $a.e. \; x\in
\Omega$, Fatou's Lemma gives
 \[
  \int_\Omega (u^\ast_i -I)^2 p^\ast_i
  \le \liminf_{n\to\infty} \int_\Omega (u^n_i-I)^2 p^n_i.
 \]

Combining both cases just analyzed above, we have established that
 \beginEL{E4:lowerG}
    \lambda \sumK \int_\Omega (u^\ast_i -I)^2 p^\ast_i +
    \alpha \sumK \int_\Omega |\grad u^\ast_i|^2
 \le \liminf_{n \to \infty}
    \lambda \sumK \int_\Omega (u^n_i -I)^2 p^n_i +
    \alpha \sumK \int_\Omega |\grad u^n_i|^2.
 \closeE
Together with (\ref{E4:lowerA}) and (\ref{E4:lowerC}), this implies
 \[
        E[\P^\ast, \U^\ast \mid I]
    \le \liminf_{n \to \infty} E[\P^n, \U^n \mid I]
     = \inf_{(\P,\U)} E[\P,\U \mid I],
 \]
and hence $(\P^\ast, \U^\ast)$ must be a minimizer. \byeproof
\end{proof}

The proof reveals an important behavior of the model~(\ref{E4:SMSa}). If certain channel
$i$ becomes dumb (i.e., $p^\ast_i \equiv 0$), it has often been introduced unnecessarily
in the first place, and the associated optimal pattern $u^\ast_i$ could be any
featureless constant image.

\begin{theorem} \label{T4:supervision}
Suppose $I \in L^2(\Omega)$. Then an optimal pattern-ownership pair must exist to the
soft Mumford-Shah segmentation model~(\ref{E4:SMSa}) with
supervision~(\ref{E4:supervision}), assuming that each patch $Q_i$ has a positive
Lebesgue measure $|Q_i|>0$.
\end{theorem}

\begin{proof}
The proof is almost identical to the unsupervised case above, and simplifies
substantially by noticing that no channel could become dumb due to supervision.
Furthermore, the  functions $\rho_i$'s in the above proof can be directly set to be
 \[ \rho_i = \frac 1 {|Q_i|} 1_{Q_i}(x), \qquad i=1:K, \]
without the necessity of turning to Lemma~\ref{L4:meanbound}.
 \byeproof
\end{proof}

\subsection{Mixture of Homogeneous Gaussians}

When each pattern $i$ is a homogeneous Gaussian $N(m_i, \sigma)$ with a distinct mean
value $m_i$, one has
 \[ u_i(x) \equiv m_i, \qquad x \in \Omega, \quad i=1:K. \]
Define $\bm=(m_1, \dots, m_K)$. As a result, the soft Mumford-Shah model~(\ref{E4:SMSa})
simplifies to
 \beginEL{E4:pcSMS}
 \min_{(\P, \bm)} E[\P, \bm \mid I]
  = \lambda \sumK \int_\Omega (I -m_i)^2 p_i +
    \sumK \int_\Omega \left( 9\eps|\grad p_i|^2 + {(p_i(1-p_i))^2 \over \eps} \right).
 \closeE

\begin{theorem}
Suppose $I\in L^2(\Omega)$. Then a minimizer pair $(\P^\ast, \bm^\ast)$ to $E[\P, \bm
\mid I]$ exists for both the unsupervised and supervised cases.
\end{theorem}

The proof can be derived readily from the previous general cases and is hence left out.
When $K=2$, a similar model was proposed earlier by Shen~\cite{shen_Gampcms} under the
symmetrization transform:
 \[
 p_1(x) = {1- z(x) \over 2 }, \quad
 p_2(x) = {1+ z(x) \over 2 }, \quad
 z \in [-1,  1].
 \]

The model~(\ref{E4:pcSMS}) could be considered as the soft version of Chan and Vese's
model~\cite{chaves_ac} from the point of view of region-based active contours. Chan and
Vese have demonstrated that such a piecewise constant Mumford-Shah model (or the CV model
as popularly referred to in the present literature) is already powerful enough for a
number of applications including medical imaging.


\section{Euler-Lagrange Equations and Computation on $(K-1)$-Simplex}

 \subsection{Euler-Lagrange Equations on $(K-1)$-Simplex}

To minimize the energy for the soft Mumford-Shah segmentation
 \beginEL{E5:SMSb}
    E[\P, \U \mid I ]
 =   \lambda \sumK \int_\Omega (u_i -I)^2 p_i
   + \alpha \sumK \int_\Omega |\grad u_i|^2
   + \sumK \int_\Omega \left( 9\eps|\grad p_i|^2 + {(p_i(1-p_i))^2 \over \eps} \right)
    ,
 \closeE
one resorts to its gradient-descent flow or Euler-Lagrange equations. In this section, we
discuss these equations and their practical computational schemes.

The first-order partial variation on $\U$ given $\P$ leads to, for $i=1:K$,
 \[
 \alpha \Delta u_i + \lambda (I-u_i)p_i=0, \quad \mbox{on $\Omega$}; \qquad
 {\pd u_i \over \pd \bn }=0, \quad \mbox{along $\pd \Omega$,}
 \]
where $\bn$ stands for the outer normal vector field along $\pd \Omega$. Thus the
Euler-Lagrange equations on the patterns are all in the form of linear Poisson equations
with variable coefficient fields:
 \[
 -\alpha \Delta u_i + (\lambda p_i) u_i =  f_i, \qquad i=1:K,
 \]
with Neumann adiabatic boundary conditions, where the source terms are $f_i(x)=\lambda
p_i(x) I(x)$.

The first-order variation on the ownerships $\P$ is carried out on the probability
$(K-1)$-simplex $\Delta_{K-1}$, which is a compact manifold (with border) of codimension
1 embedded in $\bbR^K$. Chan and Shen~\cite{chashe_nonflat} developed a general framework
for modelling and computing image features that ``live" on general manifolds, and
especially those that are embedded in $\bbR^K$. We shall follow the approach there.

Without the simplex constraint on the ownerships, for any given $\U$, the first order
variation of the soft energy $E$ under $\P \to \P +\delta \P$ is given by
 \[
 \delta E = \int_\Omega \sumK V_i \delta p_i dx + \int_{\pd \Omega} \sumK v_i \delta p_i
 d\H^1,
 \]
where $\H^1$ is the 1-D Hausdorff measure along $\pd \Omega$, and
 \begin{align}
 V_i &= \lambda (u_i -I)^2 - 18 \eps \Delta p_i + 2 \eps^{-1} p_i (1-p_i) (1-2p_i),
    \label{E5:Vflux} \\
 v_i &= 18 \eps \frac{\pd p_i}{\pd \bn}, \quad \mbox{along $\pd \Omega$}.
    \label{E5:vflux}
 \end{align}
Define $\V=(V_1, \dots, V_K)$ and $\bv=(v_1, \dots, v_K)$. Then
 \[ \delta E = \int_\Omega \V \cdot \delta\P dx + \int_{\pd \Omega} \bv \cdot \delta \P d\H^1, \]
which holds for any {\em free} variation of $\P$ in $\bbR^K$, or one writes in the
free-gradient form
 \[ {{\pd E} \over {\pd_f \P}} = \V \big|_\Omega + \bv \big|_{\pd \Omega}. \]

In reality, $\P \in \Delta_{K-1}$. Let $T_\P \Delta_{K-1}$ denote the tangent space of
$\Delta_{K-1}$ at any single point $\P \in \Delta_K$, and
 \[
 \pi: T_\P \bbR^K \to T_\P \Delta_{K-1}
 \]
the orthogonal projection onto the tangent space in $\bbR^K$. Since the normal direction
of the tangent plane is given by $\bone_K/\sqrt{K}=(1, \dots, 1)/\sqrt{K}$, the
projection operator is explicitly given by, for any $\bw \in T_\P \bbR^K$,
 \[
 \pi(\bw) = \bw - \bone_K \inner{\bw}{\bone_K}/K = \bw - \mean{\bw} \bone_K, \qquad
 \mbox{with} \quad
 \mean{\bw} =\frac 1 K \sumK w_i.
 \]
The constrained gradient of $E$ on $\Delta_{K-1}$ is therefore given by
 \[
 {\pd E \over \pd \P}
    = \pi\left( {\pd E \over \pd_f \P}  \right)
    =       (\V - \mean{\V} \bone_K) \big|_\Omega
        +   (\bv - \mean{\bv} \bone_K) \big|_{\pd \Omega}.
 \]
In particular, the system of Euler-Lagrange equations on $\P$ given $\U$ is given by
 \beginEL{E5:pEL}
 \begin{cases}
 V_i(x) - \mean{\V}(x)=0, & x \in \Omega, \\
 v_i(z) - \mean{\bv}(z)=0, & z \in \pd \Omega,
 \end{cases}
 \closeE
for $i=1:K$. The coupling among different channels is evident from these two formulae.

\begin{lemma}
Suppose $\P: \Omega \to \Delta_{K-1}$. Then for any $z\in \pd \Omega$, $\mean{\bv}(z)=0$,
where the boundary ``flux" $\bv$ is defined in~(\ref{E5:vflux}).
\end{lemma}

\begin{proof}
This is obtained by direct computation: at any $z \in \pd \Omega$,
 \[
 \begin{split}
 \mean{\bv} &= \frac 1 K \sumK v_i
     = \frac {18 \eps} K \sumK {\pd p_i \over \pd \bn}       \\
    &= \frac {18 \eps} K \frac{\pd}{\pd \bn} \left( \sumK p_i \right)
     = \frac {18 \eps} K \frac{\pd 1}{\pd \bn}=0.
 \end{split}
 \]
\byeproof
\end{proof}

As a result, the boundary conditions in~(\ref{E5:pEL}) simplify to the ordinary Neumann
conditions $\pd p_i/\pd \bn =0, i=1:K$. Combining all the above derivations, we have
established the following theorem.

\begin{theorem}[Euler-Lagrange Equations] \label{T5:ELsys}
The system of Euler-Lagrange equations of $E[\P, \U \mid I]$ are given by
 \beginEL{E5:upEL}
 \begin{split}
 -\alpha \Delta u_i +  (\lambda p_i) u_i & = (\lambda p_i) I, \\
 -18\eps\Delta p_i  + 2 \eps^{-1} p_i(1-p_i)(1-2p_i) & = \mean{\V} -\lambda (u_i -I)^2,
 \qquad i=1:K,
 \end{split}
 \closeE
on $\Omega$, all with Neumann boundary conditions along $\pd \Omega$. Here $\V=\V(\P,
\U)$ is defined as in (\ref{E5:Vflux}). Furthermore, under
supervision~(\ref{E4:supervision}), the ownerships must satisfy the interpolation
conditions:
 \[ p_i \big|_{Q_j} = \delta_{i,j}, \qquad i, j=1:K, \]
Or equivalently, the equations on $p_i$'s in (\ref{E5:upEL}) hold on $\Omega \setminus
\left( \cup_{i=1}^K Q_i \right)$ with
 \[
 \mbox{ Neumann conditions along $\pd \Omega$, and Dirichlet
       conditions along $\cup_{i=1}^K \pd Q_i: p_i \big|_{\pd Q_j} = \delta_{i,j}$.}
 \]
\end{theorem}

Similarly, one has the following result for the piecewise constant SMS
model~(\ref{E4:pcSMS}), which carries much lower complexity compared with the full SMS
model.

\begin{proposition}[Euler-Lagrange Equations for Piecewise Constant SMS]
The Euler-Lagrange equations for $E[\P, \bm \mid I]$ in~(\ref{E4:pcSMS}) are given by
\beginEL{E5:mpEL}
 \begin{split}
  m_i = \mean{I}_{p_i} := {\int_\Omega I p_i  \over \int_\Omega p_i},   \\
 -18\eps\Delta p_i  + 2 \eps^{-1} p_i(1-p_i)(1-2p_i) & = \mean{\V} -\lambda (m_i -I)^2,
 \qquad i=1:K,
 \end{split}
 \closeE
with Neumann conditions for all the ownerships $p_i$'s along $\pd \Omega$.
\end{proposition}

\subsection{Computation of the Euler-Lagrange Equations}

Computationally, as well practiced in multivariate optimization problems, (\ref{E5:upEL})
and (\ref{E5:mpEL}) can be solved via the algorithm of {\em alternating minimization}
(AM)~\cite{eseshe,she_dejitter}. The AM algorithm is closely connected to the celebrated
EM (expectation-maximization) algorithm in statistical estimation problems with hidden
variables~\cite{JepsonBlack,EM_book_Litt}. In the current context, the ownership
distributions $p_i$'s could be treated as the hidden variables.

Like EM, the AM algorithm is progressive. Given the current ($t=n$) best estimation of
the patterns $\U^n=(u^n_i \mid i=1:K)$, by solving
 \[
 \P^n =\argmin_\P \; E[\P \mid \U^n, I],
 \]
or equivalently,
 \beginEL{E5:U2P}
 -18\eps \Delta p_i + 2 \eps^{-1}p_i(1-p_i)(1-2p_i)=\mean{\V^n}-\lambda (u^n_i -I)^2, \quad
 i=1:K,
 \closeE
with Neumann boundary conditions, one obtains the current best estimation of the
ownerships $\P^n =(p^n_i \mid i=1:K)$. Subsequently, based on $\P^n$, by solving
 \[
 \U^{n+1} =\argmin E[\U \mid \P^n, I],
 \]
or equivalently,
 \beginEL{E5:P2U}
 -\alpha \Delta u_i + (\lambda p^n_i) u_i = (\lambda p^n_i) I,  \qquad i=1:K,
 \closeE
with Neumann boundary conditions, one completes a single round of pattern updating $\U^n
\to \U^{n+1}$. The same procedure applies to the piecewise constant soft Mumford-Shah
equations in (\ref{E5:mpEL}).

Since the system (\ref{E5:P2U}) is linear and decoupled, the main computational
complexity resides in the integration of $(\ref{E5:U2P})$, which is coupled and nonlinear
due to the simplex constraint and the double-well potential in the energy. Define
$e_i(x)=(u_i(x) - I(x))^2$ and $\be=(e_i \mid i=1:K)$. In order to solve
 \beginEL{E5:peq}
    -18\eps \Delta p_i + 2\eps^{-1}p_i(1-p_i)(1-2p_i)=\mean{\V} - \lambda e_i
 \closeE
given $\be$ and $\V=\V(\P, \U)=\V(\P,\be)$ (see (\ref{E5:Vflux})), first notice that
 \[
 \begin{split}
  \mean{\V(\P,\be)}
 &= \frac 1 K \sumK
    \left( -18\eps \Delta p_i + e_i + 2 \eps^{-1}p_i(1-p_i)(1-2p_i) \right) \\
 &= \frac 1 K \sumK e_i +
    {2 \eps^{-1} \over K} \sumK (2p_i^3 - 3p_i^2) + {2 \eps^{-1} \over K},
 \end{split}
 \]
since $\sumK p_i =1$ and $\Delta \left(\sumK p_i\right)=0$. We also split the
double-potential force in (\ref{E5:peq}) by
 \[
 p_i (1-p_i)(1-2p_i)= p_i(1-p_i)^2 - p_i^2 (1-p_i).
 \]
In combination, the nonlinear equation (\ref{E5:peq}) can then be solved iteratively:
 \[ \dots \to \P^{\mean{j}} \to \P^{\mean{j+1}} \to \dots \]
by the following linearization procedure:
 \[
 \begin{split}
  &  -18 \eps \Delta p_i^{\mean{j+1}} + 2 \eps^{-1} p_i^{\mean{j+1}} \left(1-p_i^{\mean{j}}\right)^2
  = f^{\mean{j}}_i, \\
  & f^{\mean{j}}_i = - e_i + \mean{V(\P^{\mean{j}}, \be)}
    + 2\eps^{-1} \left( p_i^{\mean{j}} \right)^2 \left(1-p_i^{\mean{j}}\right),
 \end{split}
 \]
with Neumann adiabatic boundary conditions for all the channels $i=1:K$. This system of
linear Poisson equations can be conveniently integrated using any elliptic solvers. The
detailed {\em numerical analysis} on the convergence rates of the entire algorithm above,
however, is still an open problem and well deserves some systematic investigation.


\section{Computational Examples}

In this section, we present the computational results of the proposed soft Mumford-Shah
model. Notice that the extension of the above SMS models to color images is
straightforward by having the gray values $u_i$'s replaced by RGB vectors. (We, however,
must remind the reader that {\em perceptually} RGB may not be the most ideal
representation of colors compared with other nonlinear approaches, e.g.,
brightness-chromaticity~\cite{chashe_nonflat} and HSV~\cite{chakanshe_hsv}.)

Figures~\ref{F6:A} and \ref{F6:B} illustrate the performance of the SMS model on two
synthetic images with multiple phases. Figure~\ref{F6:A} shows a typical T-junction and
Figure~\ref{F6:B} a 3-phase image with a narrow bottleneck. Plotted in the figures are
the hard segments obtained from the SMS model via the hardening formulae~(\ref{E1:s2hA})
and (\ref{E1:s2hB}).

{
\begin{figure}[ht]
 \centering{
 \epsfig{file=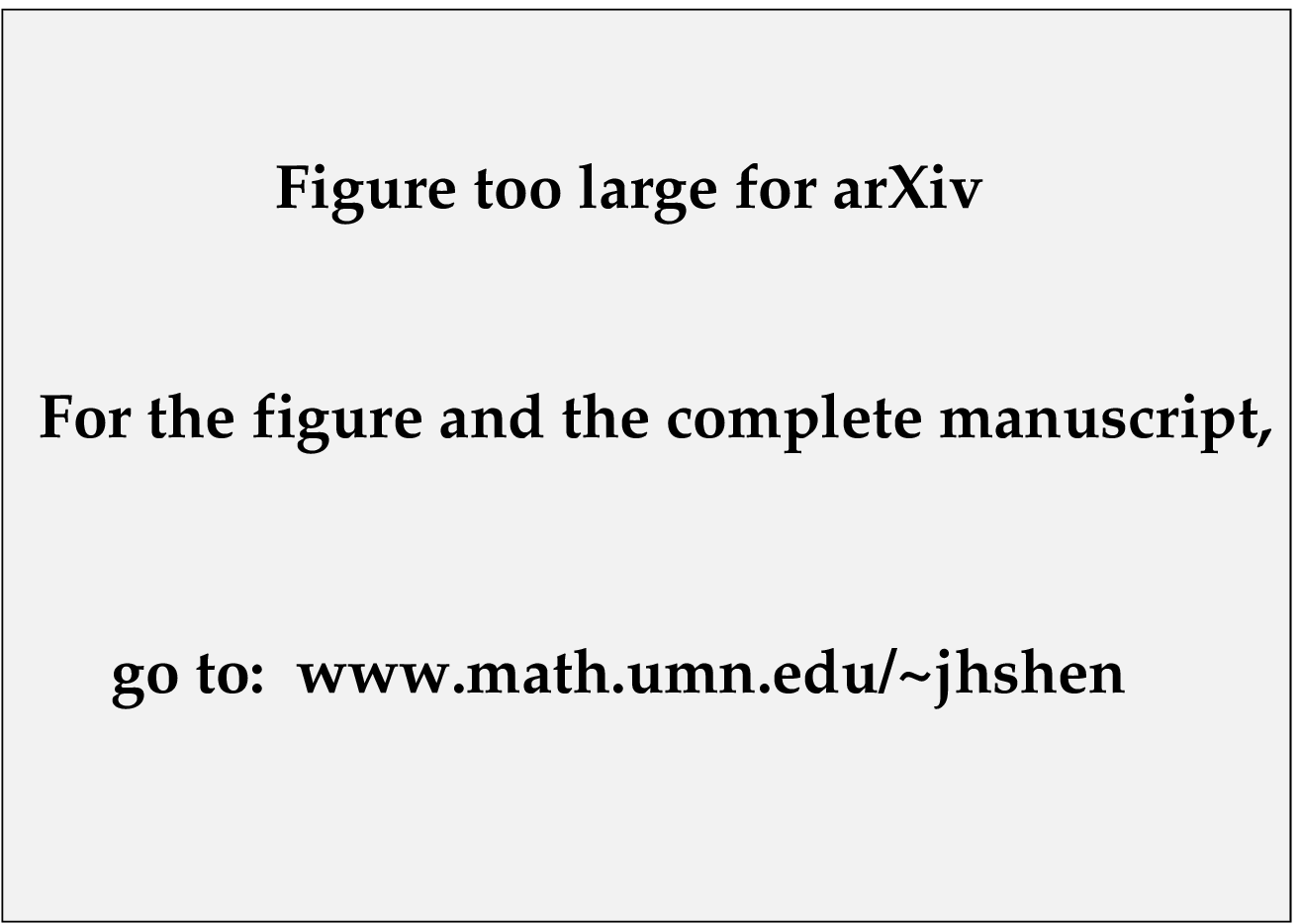, height=2.6in}
 \caption{Synthetic image of a T-junction: hard segmentation from the SMS model via ``hardening"
          formulae~(\ref{E1:s2hA}) and (\ref{E1:s2hB}). The 120-degree regularization behavior
          at the junction point is also well known in the classical Mumford-Shah model~\cite{mumsha}.}
\label{F6:A} }
\end{figure}
}

{
\begin{figure}[ht]
 \centering{
 \epsfig{file=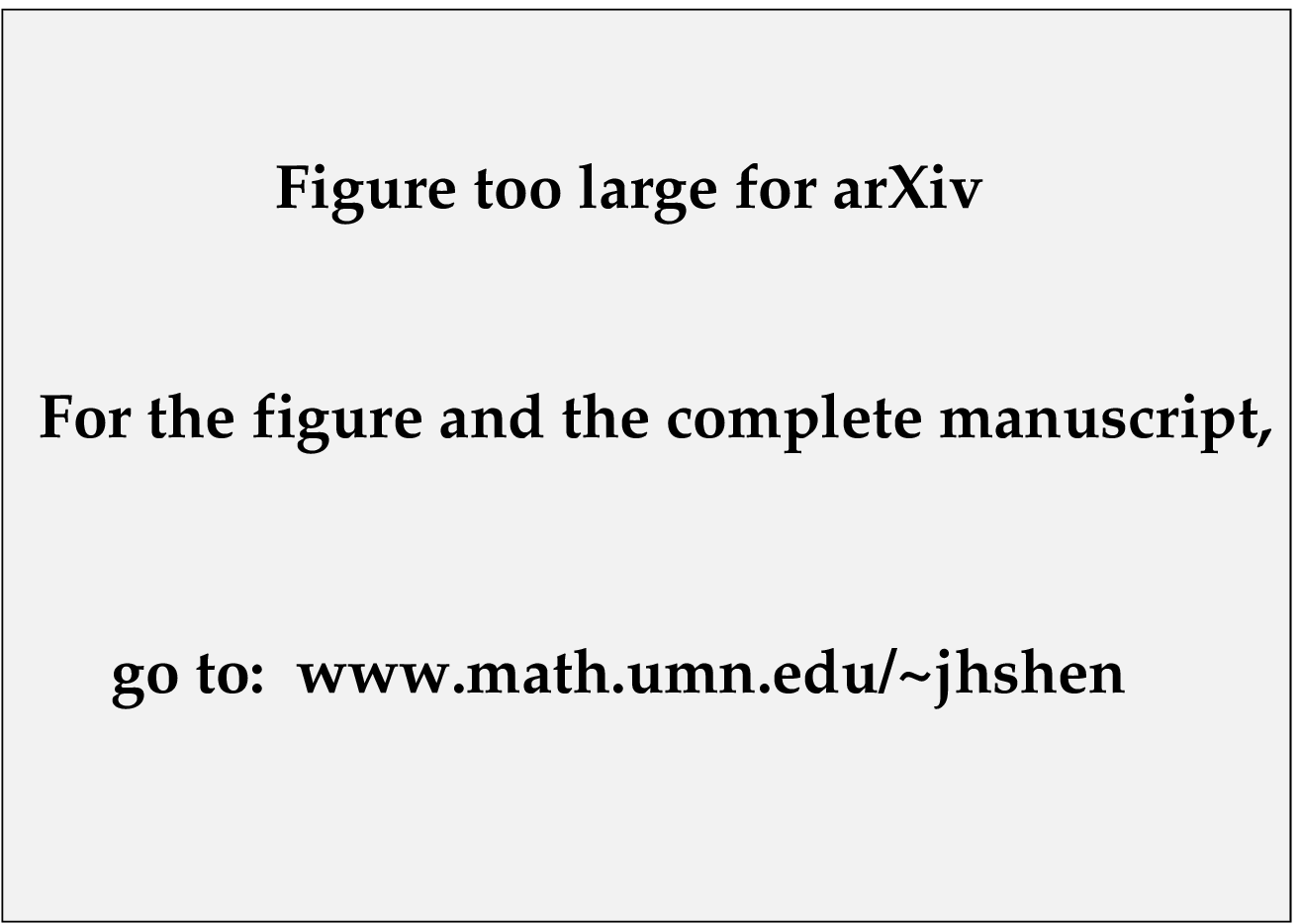, height=3in}
 \caption{Synthetic image of a narrow bottleneck: hard segmentation from the SMS model via ``hardening"
          formulae~(\ref{E1:s2hA}) and (\ref{E1:s2hB}). The thickening regularization at the
          bottleneck junction can be explained similarly by the classical Mumford-Shah model
          for which minimum-surface or ``soap-foam" behavior arises
          due to the surface tension energy. Also see the recent work by Kohn and
          Slastikov~\cite{KohnDoubleWell} for the singularity analysis of a similar
          problem arising from micromagnetism. } \label{F6:B} }
\end{figure}
}

Plotted in Figure~\ref{F6:C} are the hardened segments obtained from the soft
Mumford-Shah segmentation model via formulae~(\ref{E1:s2hA}) and (\ref{E1:s2hB}). For
this application, a user specifies three small patches (three rectangles in this example)
$Q_1, Q_2,$ and $Q_3$, and the SMS model proceeds with the extra interpolation conditions
in (\ref{E4:supervision}) for the ownerships.

{
\begin{figure}[ht]
 \centering{
 \epsfig{file=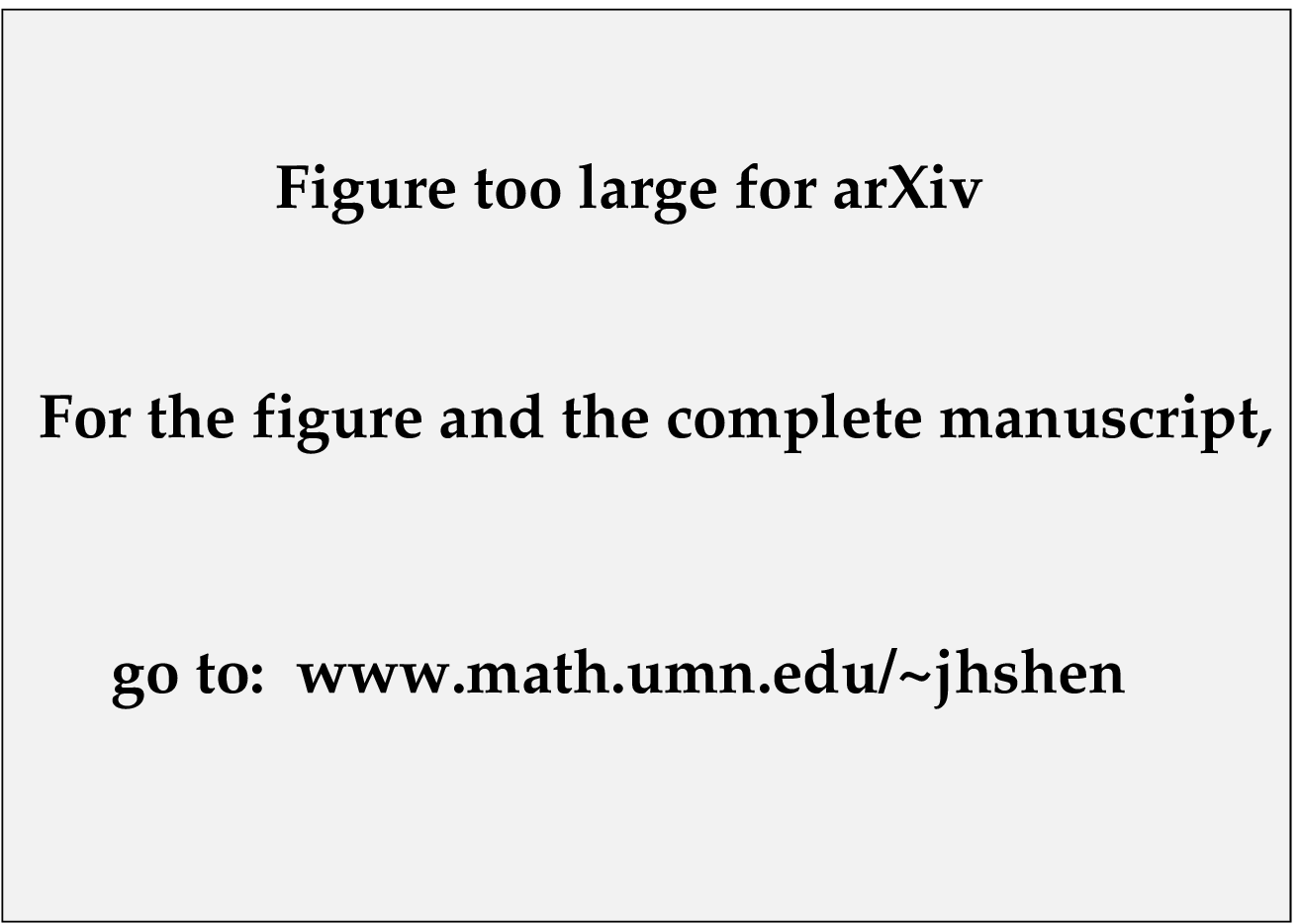, height=3in}
 \caption{A real noisy brain image: hard segmentation from the SMS model via ``hardening"
          formulae~(\ref{E1:s2hA}) and (\ref{E1:s2hB}).  }
\label{F6:C} }
\end{figure}
}


In Figure~\ref{F6:D}, another example of a natural image is segmented via the SMS model
and the ``hardening" formulae~(\ref{E1:s2hA}) and (\ref{E1:s2hB}). A user supervises with
three patches $Q_1, Q_2$ and $Q_3$, and designates the two on the body to a pattern
ownership $p_{\mathrm{body}}$ and the third (from the ocean) to $p_{\mathrm{ocean}}$. If
the three are treated as distinct patterns, the SMS model still works, but one needs an
extra step of high-level vision processing (e.g., based on Grenander's graph
models~\cite{gre}) to group the skin-tone and the purple-shirt patterns in order to
capture the entire body faithfully.

{
\begin{figure}[ht]
 \centering{
 \epsfig{file=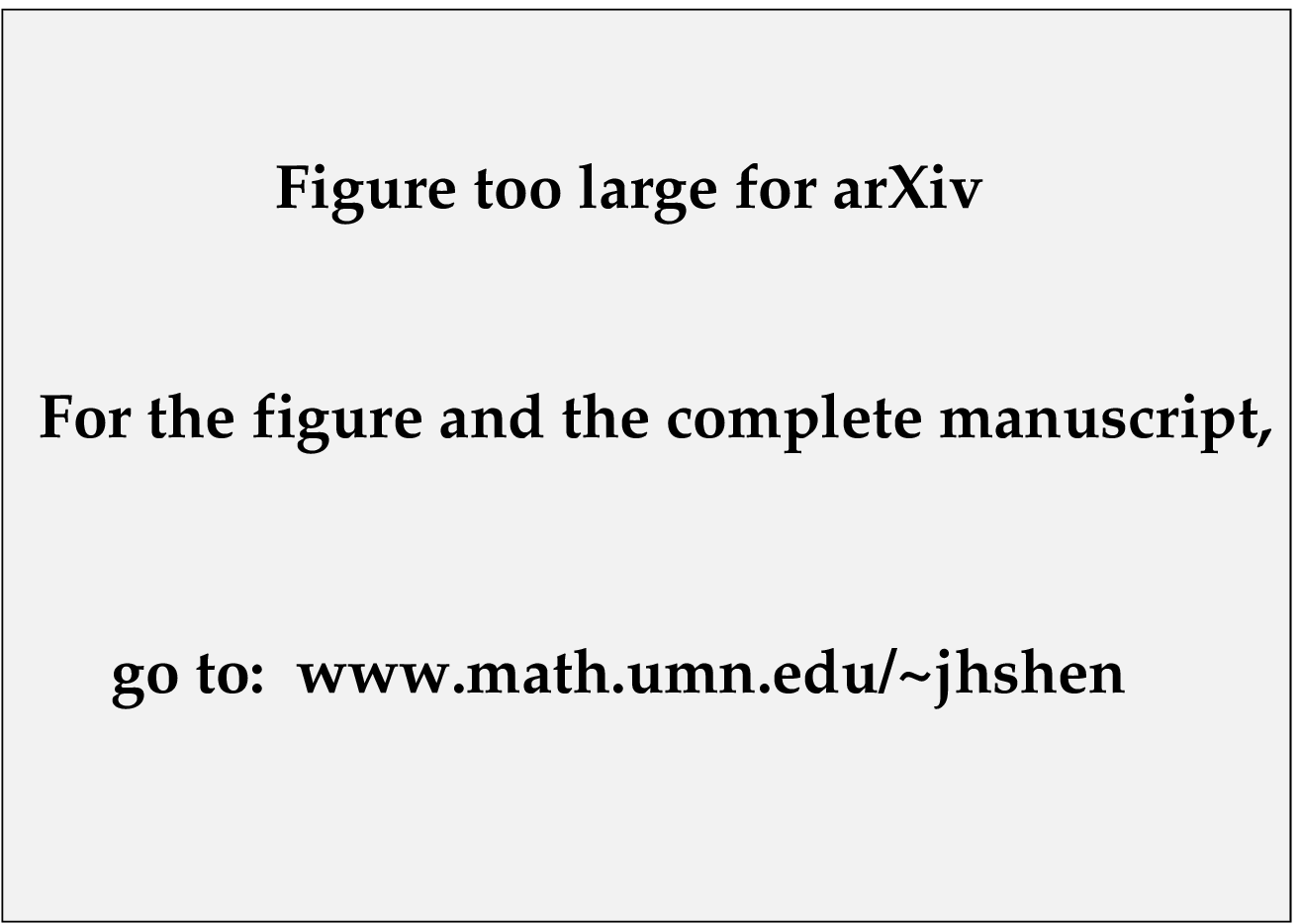, height=2in}
 \caption{ Hard segmentation from the SMS model via ``hardening"
          formulae~(\ref{E1:s2hA}) and (\ref{E1:s2hB}), based on a
          2-phase supervision. Denote the two rectangles on the body by
          $Q_1$ and $Q_2$, and the third by $Q_3$. Supervision provides
          the ownership interpolation condition: $p_{\mathrm{body}}=1$ on
          $Q_1 \cup Q_2$ and $0$ on $Q_3$, while $p_{\mathrm{ocean}}=1$ on $Q_3$
          and $0$ on $Q_1 \cup Q_2$. Patch selection can also be automated
          based on multiscale patch statistics (e.g.,  Li and Perona~\cite{feifei05}). }
 \label{F6:D} }
\end{figure}
}

Finally, plotted in Figures~(\ref{F6:E}) and (\ref{F6:F}) are the ownerships from the SMS
model based on the 3-phase and 4-phase supervision separately in
Figure~\ref{F4:superbeach}. The stochastic nature of the outcomes (i.e., the softly
transiting ownerships $p_i$'s instead of hard segmentation) is closer to the way how a
human subject may perceive such a natural scene. In particular, the SMS model seems to be
consistent with the most recent theory that {\em hard} pattern segments may not be
absolutely necessary for natural scene recognition~\cite{feifei99,feifei05}.

{
\begin{figure}[ht]
 \centering{
 \epsfig{file=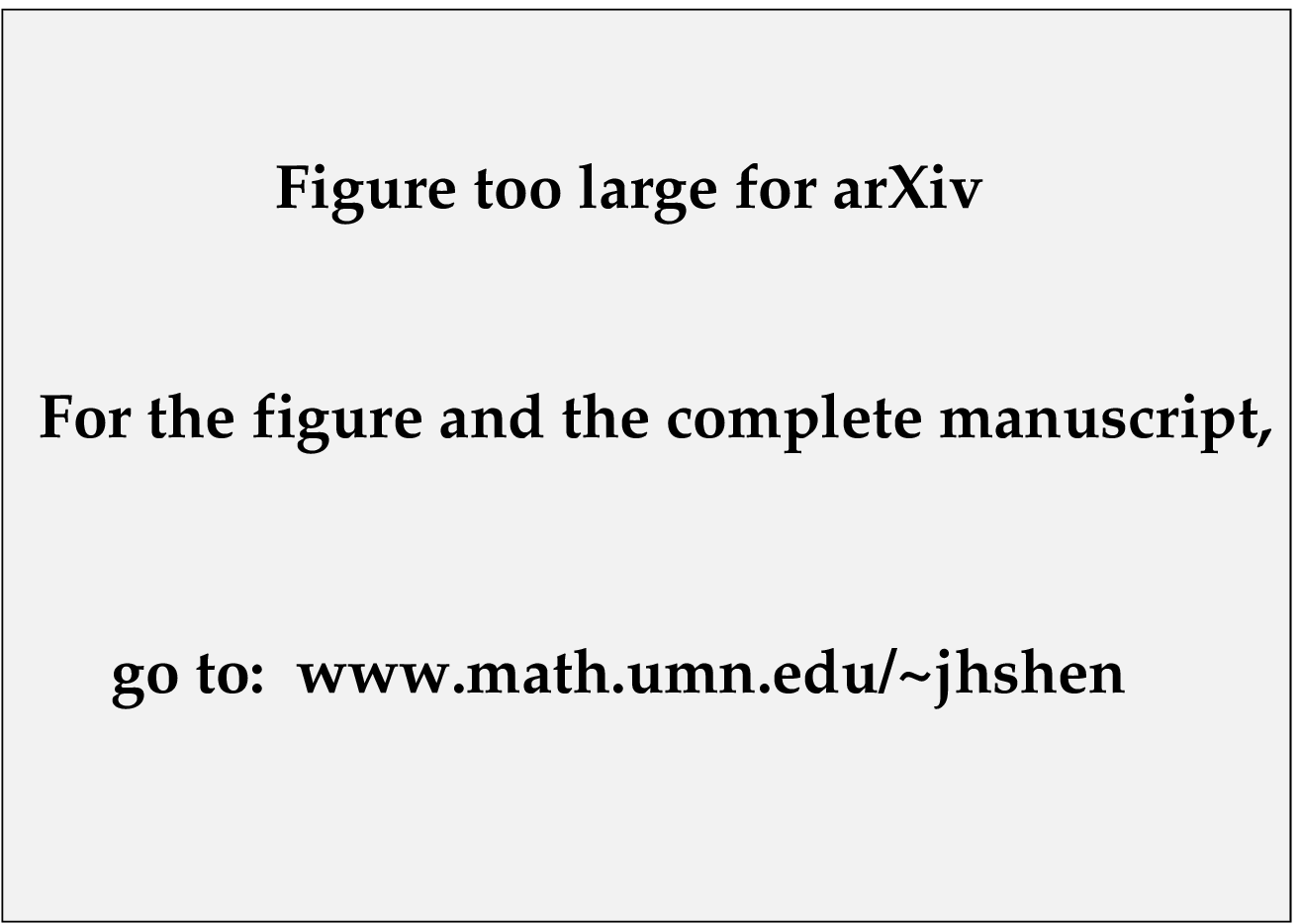, height=3in}
 \caption{Soft Mumford-Shah segmentation with three phases corresponding to
          the supervision on the {\em left} panel of Figure~\ref{F4:superbeach}. Plotted here are
          the three ownership distributions $p_1(x), p_2(x)$, and $p_3(x)$. Due to ``under"-supervision,
          namely the number $K$ of specified patterns is less than that of
          the visually meaningful ones, the grass pattern has ``absorbed"
          the ocean pattern due to the greenish color they happen to share.}
\label{F6:E} }
\end{figure}
}

{
\begin{figure}[ht]
 \centering{
 \epsfig{file=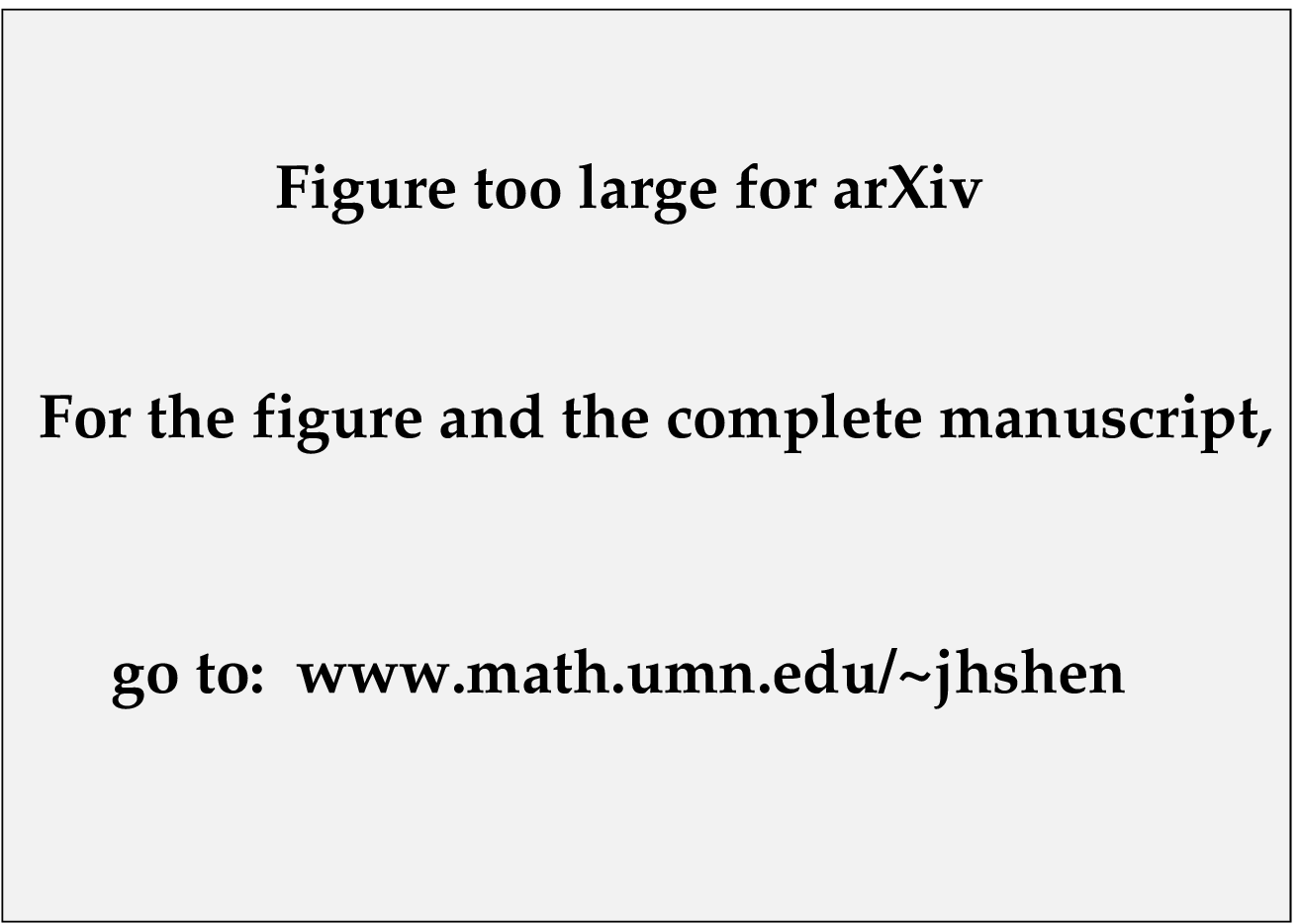, height=3in}
 \caption{Soft Mumford-Shah segmentation with four phases corresponding to
          the supervision on the {\em right} panel of Figure~\ref{F4:superbeach}.
          Plotted here are the four ownership
          distributions $p_1(x), p_2(x), p_3(x)$, and $p_4(x)$. Unlike
          Figure~\ref{F6:E}, the narrow ocean pattern is now softly
          segmented  due to the extra fourth patch $Q_4$.}
\label{F6:F} }
\end{figure}
}


\section*{Acknowledgments}

The author is very grateful to Prof. Alan Yuille for an enlightening discussion after the
current work was first presented. For their generous teaching and continual inspiration,
the author is always profoundly indebted to Professors Gil Strang, Tony Chan, Stan Osher,
David Mumford, Jean-Michel Morel, and Stu Geman. The author must thank his wonderful
former teacher, Prof. Dan Kerstan at the Psychology Department of the University of
Minnesota, for his first introduction on mixture image models and stochastic visual
processing several year ago. The author also thanks the Institute of Mathematics and its
Applications ({\bf IMA}) and the Institute of Pure and Applied Mathematics ({\bf IPAM})
for their persistent role in supporting this new emerging field. Finally, the author
would like to dedicate this paper to his dear friends Yingnian Wu and Song-Chun Zhu for
the unique friendship cultivated by the intellectually rich soil of vision and cognitive
sciences.


\end{document}